% 2001
% 20 Mars
% 27 mars : corrections
% 3 avril : corrections
% 30 mai : enlev\'e la capitale dans Arithmetischen dans le titre
% 20 juin 2001 : enlev\'e sur Rost F_4, a completer
% 1 mars 2007 
\input amssym.def 
\input amssym.tex
%Reglages generaux ÑÑÑÑÑÑÑÑÑÑÑÑÑÑÑÑÑÑÑÑÑÑÑÑÑ
% ÑÑÑÑÑÑÑÑÑÑÑÑÑÑÑÑÑÑÑÑÑ que vous pouvez modifier ÑÑÑÑÑÑÑÑÑÑÑÑÑÑÑÑÑÑ
%ÑÑÑÑÑÑÑÑÑÑÑÑÑÑÑÑÑÑÑÑÑÑ sans faire de degats ÑÑÑÑÑÑÑÑÑÑÑÑÑÑÑÑÑÑÑÑÑÑ

\hsize=13.50cm    
\vsize=18cm       
\parindent=12pt   \parskip=0pt      
\pageno=1 

% deux jeux d'offsets. Le premier concerne la \magnification=1000
% (si l'image n'est pas agrandie)

\hoffset=15mm    % offset horizontal en \magnification=1000
\voffset=1cm    % offset vertical en \magnification=1000
 
% Le deuxieme concerne  la \magnification=1200
% (lorsque l'image est agrandie de 20%) 

\ifnum\mag=\magstep1
\hoffset=-2mm   % offset horizontal en \magnification=1200
\voffset=.8cm   % offset horizontal en \magnification=1200
\fi

% ÑÑÑÑÑÑÑÑÑÑÑÑÑÑÑÑÑÑÑÑÑÑÑÑÑÑÑ Reglages ÑÑÑÑÑÑÑÑÑÑÑÑÑÑÑÑÑÑÑÑÑÑÑÑÑÑ
% ÑÑÑÑÑÑÑÑÑÑÑ auquels il ne vaut mieux pas toucher ÑÑÑÑÑÑÑÑÑÑÑÑÑÑ

\pretolerance=500 \tolerance=1000  \brokenpenalty=5000

% ÑÑÑÑÑÑÑÑÑÑÑÑÑÑÑÑÑÑÑÑ Debut des macros privees ÑÑÑÑÑÑÑÑÑÑÑÑÑÑÑÑÑÑÑ
\catcode`\@=11
% ÑÑÑÑÑÑÑÑÑÑÑÑÑÑÑÑÑÑÑÑÑ Les fontes ÑÑÑÑÑÑÑÑÑÑÑÑÑÑÑÑÑÑÑÑÑÑÑÑÑÑÑÑÑÑÑÑ

\font\eightrm=cmr8         \font\eighti=cmmi8
\font\eightsy=cmsy8        \font\eightbf=cmbx8
\font\eighttt=cmtt8        \font\eightit=cmti8
\font\eightsl=cmsl8        \font\sixrm=cmr6
\font\sixi=cmmi6           \font\sixsy=cmsy6
\font\sixbf=cmbx6

% Fontes AMS

\font\tengoth=eufm10       \font\tenbboard=msbm10
\font\eightgoth=eufm10 at 8pt      \font\eightbboard=msbm10 at 8 pt
\font\sevengoth=eufm7      \font\sevenbboard=msbm7
\font\sixgoth=eufm7 at 6 pt        \font\fivegoth=eufm5

 \font\tencyr=wncyr10       
\font\eightcyr=wncyr10 at 8 pt      
\font\sevencyr=wncyr10 at 7 pt      
\font\sixcyr=wncyr10 at 6 pt

% Pour que les accents se placent correctement en mode math en corps 8 et 6

\skewchar\eighti='177 \skewchar\sixi='177
\skewchar\eightsy='60 \skewchar\sixsy='60

% Nouvelles familles pour les maths

\newfam\gothfam           \newfam\bboardfam
\newfam\cyrfam

\def\tenpoint{%
  \textfont0=\tenrm \scriptfont0=\sevenrm \scriptscriptfont0=\fiverm
  \def\rm{\fam\z@\tenrm}%
  \textfont1=\teni  \scriptfont1=\seveni  \scriptscriptfont1=\fivei
  \def\oldstyle{\fam\@ne\teni}\let\old=\oldstyle
  \textfont2=\tensy \scriptfont2=\sevensy \scriptscriptfont2=\fivesy
  \textfont\gothfam=\tengoth \scriptfont\gothfam=\sevengoth
  \scriptscriptfont\gothfam=\fivegoth
  \def\goth{\fam\gothfam\tengoth}%
  \textfont\bboardfam=\tenbboard \scriptfont\bboardfam=\sevenbboard
  \scriptscriptfont\bboardfam=\sevenbboard
  \def\bb{\fam\bboardfam\tenbboard}%
 \textfont\cyrfam=\tencyr \scriptfont\cyrfam=\sevencyr
  \scriptscriptfont\cyrfam=\sixcyr
  \def\cyr{\fam\cyrfam\tencyr}%
  \textfont\itfam=\tenit
  \def\it{\fam\itfam\tenit}%
  \textfont\slfam=\tensl
  \def\sl{\fam\slfam\tensl}%
  \textfont\bffam=\tenbf \scriptfont\bffam=\sevenbf
  \scriptscriptfont\bffam=\fivebf
  \def\bf{\fam\bffam\tenbf}%
  \textfont\ttfam=\tentt
  \def\tt{\fam\ttfam\tentt}%
  \abovedisplayskip=12pt plus 3pt minus 9pt
  \belowdisplayskip=\abovedisplayskip
  \abovedisplayshortskip=0pt plus 3pt
  \belowdisplayshortskip=4pt plus 3pt 
  \smallskipamount=3pt plus 1pt minus 1pt
  \medskipamount=6pt plus 2pt minus 2pt
  \bigskipamount=12pt plus 4pt minus 4pt
  \normalbaselineskip=12pt
  \setbox\strutbox=\hbox{\vrule height8.5pt depth3.5pt width0pt}%
  \let\bigf@nt=\tenrm       \let\smallf@nt=\sevenrm
  \normalbaselines\rm}

\def\eightpoint{%
  \textfont0=\eightrm \scriptfont0=\sixrm \scriptscriptfont0=\fiverm
  \def\rm{\fam\z@\eightrm}%
  \textfont1=\eighti  \scriptfont1=\sixi  \scriptscriptfont1=\fivei
  \def\oldstyle{\fam\@ne\eighti}\let\old=\oldstyle
  \textfont2=\eightsy \scriptfont2=\sixsy \scriptscriptfont2=\fivesy
  \textfont\gothfam=\eightgoth \scriptfont\gothfam=\sixgoth
  \scriptscriptfont\gothfam=\fivegoth
  \def\goth{\fam\gothfam\eightgoth}%
  \textfont\cyrfam=\eightcyr \scriptfont\cyrfam=\sixcyr
  \scriptscriptfont\cyrfam=\sixcyr
  \def\cyr{\fam\cyrfam\eightcyr}%
  \textfont\bboardfam=\eightbboard \scriptfont\bboardfam=\sevenbboard
  \scriptscriptfont\bboardfam=\sevenbboard
  \def\bb{\fam\bboardfam}%
  \textfont\itfam=\eightit
  \def\it{\fam\itfam\eightit}%
  \textfont\slfam=\eightsl
  \def\sl{\fam\slfam\eightsl}%
  \textfont\bffam=\eightbf \scriptfont\bffam=\sixbf
  \scriptscriptfont\bffam=\fivebf
  \def\bf{\fam\bffam\eightbf}%
  \textfont\ttfam=\eighttt
  \def\tt{\fam\ttfam\eighttt}%
  \abovedisplayskip=9pt plus 3pt minus 9pt
  \belowdisplayskip=\abovedisplayskip
  \abovedisplayshortskip=0pt plus 3pt
  \belowdisplayshortskip=3pt plus 3pt 
  \smallskipamount=2pt plus 1pt minus 1pt
  \medskipamount=4pt plus 2pt minus 1pt
  \bigskipamount=9pt plus 3pt minus 3pt
  \normalbaselineskip=9pt
  \setbox\strutbox=\hbox{\vrule height7pt depth2pt width0pt}%
  \let\bigf@nt=\eightrm     \let\smallf@nt=\sixrm
  \normalbaselines\rm}

\tenpoint

% Definition des petites capitales qui reagissent au \tenpoint et \eightpoint
% La syntaxe est celle d'un changement de fonte :
% {\pc FERMAT}, {\pc EUCLIDE} et {\pc G\"ODEL}.

\def\pc#1{\bigf@nt#1\smallf@nt}         \def\pd#1 {{\pc#1} }

%ÑÑÑÑÑÑÑÑÑÑÑÑÑÑÑÑÑÑÑÑÑÑ dactylographie francaise ÑÑÑÑÑÑÑÑÑÑÑÑÑÑÑÑÑ

\catcode`\;=\active
\def;{\relax\ifhmode\ifdim\lastskip>\z@\unskip\fi
\kern\fontdimen2  -1.2 \fontdimen3 \string;}

\catcode`\:=\active
\def:{\relax\ifhmode\ifdim\lastskip>\z@\unskip\fi\penalty\@M\ \fi\string:}

\catcode`\!=\active
\def!{\relax\ifhmode\ifdim\lastskip>\z@
\unskip\fi\kern\fontdimen2  -1.1 \fontdimen3 \string!}

\catcode`\?=\active
\def?{\relax\ifhmode\ifdim\lastskip>\z@
\unskip\fi\kern\fontdimen2  -1.1 \fontdimen3 \string?}

\def\^#1{\if#1i{\accent"5E\i}\else{\accent"5E #1}\fi}
\def\"#1{\if#1i{\accent"7F\i}\else{\accent"7F #1}\fi}

\frenchspacing

% ÑÑÑÑÑÑÑÑÑÑÑÑÑÑÑÑ Le format de sortie ÑÑÑÑÑÑÑÑÑÑÑÑÑÑÑÑÑÑÑÑÑÑÑÑÑÑÑÑ
% Haut et bas de page 

\newtoks\auteurcourant      \auteurcourant={\hfil}
\newtoks\titrecourant       \titrecourant={\hfil}

\newtoks\hautpagetitre      \hautpagetitre={\hfil}
\newtoks\baspagetitre       \baspagetitre={\hfil}

\newtoks\hautpagegauche     
\hautpagegauche={\eightpoint\rlap{\folio}\hfil\the\auteurcourant\hfil}
\newtoks\hautpagedroite     
\hautpagedroite={\eightpoint\hfil\the\titrecourant\hfil\llap{\folio}}

\newtoks\baspagegauche      \baspagegauche={\hfil} 
\newtoks\baspagedroite      \baspagedroite={\hfil}

\newif\ifpagetitre          \pagetitretrue  

% \nopagenumbers : c'est un peu violent, mais a marche. Alors ...

\headline={\ifpagetitre\the\hautpagetitre
\else\ifodd\pageno\the\hautpagedroite\else\the\hautpagegauche\fi\fi}

\footline={\ifpagetitre\the\baspagetitre\else
\ifodd\pageno\the\baspagedroite\else\the\baspagegauche\fi\fi
\global\pagetitrefalse}

% Redefinition de \raggedbottom pour avoir plus de mou en bas de page
% (necesssaire quand il y a beaucoup de grumeaux, des grosses
% formules centrees et pas beaucoup de texte entre)

\def\raggedbottom{\topskip 10pt plus 36pt\r@ggedbottomtrue}

% ÑÑÑÑÑÑÑÑÑÑÑÑÑÑÑÑÑÑ Macros de mise en page ÑÑÑÑÑÑÑÑÑÑÑÑÑÑÑÑÑÑÑÑÑ

% Un point-tiret

\def\pointir{\unskip . --- \ignorespaces}

% Macros Bigbreak et \Medbreak pour que les blancs verticaux ne s'ajoutent pas

\def\Bigbreak{\vskip-\lastskip\bigbreak}
\def\Medbreak{\vskip-\lastskip\medbreak}

% Texte centre dans une boite : le resultat est le plus petit rectangle
% qui contient le texte. Ce resultat est
% place dans une boite centree.
% La syntaxe est celle d'un tableau \`a une colonne

\def\ctexte#1\endctexte{%
  \hbox{$\vcenter{\halign{\hfill##\hfill\crcr#1\crcr}}$}}

% Titres centres (en gras)

\long\def\ctitre#1\endctitre{%
    \ifdim\lastskip<24pt\vskip-\lastskip\bigbreak\bigbreak\fi
  		\vbox{\parindent=0pt\leftskip=0pt plus 1fill
          \rightskip=\leftskip
          \parfillskip=0pt\bf#1\par}
    \bigskip\nobreak}

\long\def\section#1\endsection{%
\vskip 0pt plus 3\normalbaselineskip
\penalty-250
\vskip 0pt plus -3\normalbaselineskip
\Bigbreak
\message{[section \string: #1]}{\bf#1\unskip}\pointir}

\long\def\sectiona#1\endsection{%
\vskip 0pt plus 3\normalbaselineskip
\penalty-250
\vskip 0pt plus -3\normalbaselineskip
\Bigbreak
\message{[sectiona \string: #1]}%
{\bf#1}\medskip\nobreak}

\long\def\subsection#1\endsubsection{%
\Medbreak
{\it#1\unskip}\pointir}

\long\def\subsectiona#1\endsubsection{%
\Medbreak
{\it#1}\par\nobreak}

\def\rem#1\endrem{%
\Medbreak
{\it#1\unskip} : }

\def\remp#1\endrem{%
\Medbreak
{\pc #1\unskip}\pointir}

\def\rema#1\endrem{%
\Medbreak
{\it #1}\par\nobreak}

\def\newparwithcolon#1\endnewparwithcolon{
\Medbreak
{#1\unskip} : }

\def\newparwithpointir#1\endnewparwithpointir{
\Medbreak
{#1\unskip}\pointir}

\def\newpara#1\endnewpar{
\Medbreak
{#1\unskip}\smallskip\nobreak}

%ÑÑÑÑÑÑÑÑÑÑÑÑÑÑÑÑÑÑÑÑÑÑÑÑÑÑÑÑÑÑÑÑÑÑÑÑÑÑÑÑÑÑÑÑÑ
% enonces de theoremes avec numerotation apres 
% #1 = THEOREME, COROLLAIRE, etc.
% #2 = numero (par exemple 3, 3.1, etc.)
% #3 = l'enonce du th proprememnt dit.

\long\def\th#1 #2\enonce#3\endth{%
   \Medbreak
   {\pc#1} {#2\unskip}\pointir{\it #3}\medskip}

\long\def\tha#1 #2\enonce#3\endth{%
   \Medbreak
   {\pc#1} {#2\unskip}\par\nobreak{\it #3}\medskip}

%ÑÑÑÑÑÑÑÑÑÑÑÑÑÑÑÑÑÑÑÑÑÑÑÑÑÑÑÑÑÑÑÑÑÑÑÑÑÑÑÑÑÑÑÑÑ
% enonces de theoremes avec numerotation d'abord
% #1 = numero (par exemple 3, 3.1, etc.)
% #2 = THEOREME, COROLLAIRE, etc.
% #3 = l'enonce du th proprememnt dit.

\long\def\Th#1 #2 #3\enonce#4\endth{%
   \Medbreak
   #1 {\pc#2} {#3\unskip}\pointir{\it #4}\medskip}

\long\def\Tha#1 #2 #3\enonce#4\endth{%
   \Medbreak
   #1 {\pc#2} #3\par\nobreak{\it #4}\medskip}
%ÑÑÑÑÑÑÑÑÑÑÑÑÑÑÑÑÑÑÑÑÑÑÑÑÑÑÑÑÑÑÑÑÑÑÑÑÑÑÑÑÑÑÑÑÑ

% les differents retraits, voir aussi \item

\def\decale#1{\smallbreak\hskip 28pt\llap{#1}\kern 5pt}
\def\decaledecale#1{\smallbreak\hskip 34pt\llap{#1}\kern 5pt}
\def\puce{\smallbreak\hskip 6pt{$\scriptstyle\bullet$}\kern 5pt}

% ÑÑÑÑÑÑÑÑÑÑÑÑÑÑÑÑÑÑÑÑÑÑÑÑÑÑÑÑÑÑÑÑÑÑÑÑÑÑÑÑÑÑÑÑÑÑÑÑÑÑÑÑÑÑÑÑÑÑÑÑÑÑÑÑÑ
% ÑÑÑÑÑÑÑÑÑÑÑÑÑÑ ce que Knuth n'a pas fait ÑÑÑÑÑÑÑÑÑÑÑÑÑÑÑÑÑÑÑÑÑÑÑÑ
% ÑÑÑÑÑÑÑÑÑÑÑÑÑÑÑÑÑÑÑÑÑÑÑÑÑÑÑÑÑÑÑÑÑÑÑÑÑÑÑÑÑÑÑÑÑÑÑÑÑÑÑÑÑÑÑÑÑÑÑÑÑÑÑÑÑ

% Un \displaylines qui numerote ˆ droite.
% La syntaxe est la meme que celle de \eqalignno

\def\displaylinesno#1{\displ@y\halign{
\hbox to\displaywidth{$\@lign\hfil\displaystyle##\hfil$}&
\llap{$##$}\crcr#1\crcr}}

% Un \displaylines qui numerote a gauche.
% La syntaxe est la meme que celle de \leqalignno

\def\ldisplaylinesno#1{\displ@y\halign{ 
\hbox to\displaywidth{$\@lign\hfil\displaystyle##\hfil$}&
\kern-\displaywidth\rlap{$##$}\tabskip\displaywidth\crcr#1\crcr}}

% Un \eqalign qui accepte plusieurs alignements verticaux
% Motif : \hfil ** & ** \hfil & \hfil ** & ** \hfill, etc.

\def\eqalign#1{\null\,\vcenter{\openup\jot\m@th\ialign{
\strut\hfil$\displaystyle{##}$&$\displaystyle{{}##}$\hfil
&&\quad\strut\hfil$\displaystyle{##}$&$\displaystyle{{}##}$\hfil
\crcr#1\crcr}}\,}

% Systeme d'equations precede d'une accolade.
% Copi\'e sur \eqalign, on s'en sert comme une matrice
% syntaxe : signe & coef & inconnue
% les coef sont justifi\'es \`a droite (\hfil coef)
% et les inconnues \`a gauche (inconnue\hfil)
% attention : un seul & ou deux && avant le signe =
% selon la justification choisie !
% Exemple : $$\system{
%             &2 &x &- &3 & y & = &&  -5 \cr
%            -&  &x &+ &  & y & = &&   6 \cr
%           }$$

\def\system#1{\left\{\null\,\vcenter{\openup1\jot\m@th
\ialign{\strut$##$&\hfil$##$&$##$\hfil&&
        \enskip$##$\enskip&\hfil$##$&$##$\hfil\crcr#1\crcr}}\right.}

% pour avoir des messages raisonnables avec les lettres accentu\'ees

\let\@ldmessage=\message

\def\message#1{{\def\pc{\string\pc\space}%
                \def\'{\string'}\def\`{\string`}%
                \def\^{\string^}\def\"{\string"}%
                \@ldmessage{#1}}}

% ÑÑÑÑÑÑÑÑÑÑÑÑÑÑÑÑÑ Divers gadgets ÑÑÑÑÑÑÑÑÑÑÑÑÑÑÑÑÑÑÑÑÑÑÑÑÑÑÑÑÑÑÑÑ

% Pour se rendre la vie facile : \up{er}, \up{i\`eme}, n\up{0}, etc.

\def\up#1{\raise 1ex\hbox{\smallf@nt#1}}

% Utilisation : \cf. \etc.

\def\qed{\raise -2pt\hbox{\vrule\vbox to 10pt{\hrule width 4pt
                 \vfill\hrule}\vrule}}

\def\virg{\raise .4ex\hbox{,}}   % virgule aprs une fraction

 % point-virgule de ponctuation en maths

\def\build#1_#2^#3{\mathrel{
\mathop{\kern 0pt#1}\limits_{#2}^{#3}}}

% Entoure #2 d'un filet. Le filet est ecarte tout autour de #1
% Syntaxe \boxit{5pt}{...}. La ligne de base n'est pas perdue.

\def\boxit#1#2{%
\setbox1=\hbox{\kern#1{#2}\kern#1}%
\dimen1=\ht1 \advance\dimen1 by #1 \dimen2=\dp1 \advance\dimen2 by #1 
\setbox1=\hbox{\vrule height\dimen1 depth\dimen2\box1\vrule}%  
\setbox1=\vbox{\hrule\box1\hrule}%
\advance\dimen1 by .6pt \ht1=\dimen1 
\advance\dimen2 by .6pt \dp1=\dimen2  \box1\relax}

% ÑÑÑÑÑÑÑÑÑÑÑÑÑÑÑÑÑÑÑÑÑÑÑÑÑÑÑÑÑÑÑÑÑÑÑÑÑÑÑÑÑÑÑÑÑÑÑÑÑÑÑÑÑÑÑÑÑÑÑÑÑÑÑÑÑ
% fin des macros privees
% ÑÑÑÑÑÑÑÑÑÑÑÑÑÑÑÑÑÑÑÑÑÑÑÑÑÑÑÑÑÑÑÑÑÑÑÑÑÑÑÑÑÑÑÑÑÑÑÑÑÑÑÑÑÑÑÑÑÑÑÑÑÑÑÑÑ

\catcode`\@=12

% pour qu'il la ferme
\showboxbreadth=-1  \showboxdepth=-1

% \setbox1=\ctexte
% \bf NOUVEAU FORMAT a4\cr
% \noalign{\medskip}
% \bf \date\cr
% \endctexte

% $$\boxit{7pt}{\boxit{7pt}{\box1}}$$

\overfullrule=0pt

\magnification=\magstep1
\hsize=17,5truecm 
 \vsize=25.5truecm 
\hoffset=-0.9truecm 
\voffset=-0.8truecm
%\voffset=-0.3truecm
%\nopagenumbers    
%\pagenumbers
\topskip=1truecm
\footline={\tenrm\hfil\folio\hfil}
\raggedbottom
\abovedisplayskip=3mm %Reduction of space between text and formulae
\belowdisplayskip=3mm
\abovedisplayshortskip=0mm
\belowdisplayshortskip=2mm
\normalbaselineskip=12pt  %This is default. Do NOT change it!!
\normalbaselines

\def\dim{{\rm dim}}

\def\N{{\bf N}}
\def\R{{\bf R}}
\def\C{{\bf C}}
\def\F{{\bf F}}
\def\Q{{\bf Q}}
\def\Z{{\bf Z}}
\def\br{{\rm Br}}
\def\pic{{\rm Pic}}
\def\P{{\bf P}}
\def\A{{\bf A}}
\def\G{{\bf G}}
\def\gal{  {\rm Gal}   }
\def\g{{\cal G}}
\def\et{\hbox{\sevenrm \'et}   }
\def\ind{{\rm Index}}
\def\exp{{\rm Exp}}
\def\pic{{\rm Pic}}

\bigskip

Stuttgarter Tagung zum 100. Jahrestag  Richard Brauers

22.-24. M\"arz 2001

\vskip5cm

{\it Die Brauersche Gruppe; ihre Verallgemeinerungen
und  Anwendungen in der arithmetischen Geometrie}

\vskip5cm

J.-L. Colliot-Th\'el\`ene

\bigskip

C.N.R.S.

Universit\'e de Paris-Sud

Orsay

France

\vfill\eject

{\bf I. Einige Urbegriffe und Urs\"atze}
\bigskip
 (Ur = vor Grothendieck)

%Brauer, Noether, Hasse, Albert, Witt, Teichm\"uller 

\medskip

%Ch\^atelet, Roquette, Amitsur, Faddeev

Sei $k$ ein K\"orper, $k_s$ ein separabler Abschlu{\ss},
$\g=\gal(k_s/k)$.

 Unter dem Tensorprodukt bilden  die zentralen einfachen Algebren eine
abelsche Torsionsgruppe,
 die {\it  Brauersche Gruppe} 
oder {\it Brauergruppe}  von $k$, geschrieben $\br(k)$.

\bigskip

Der  Standpunkt der Faktorensystemen (Brauer, Noether) (der eng mit den
verschr\"ankten Produkten verkn\"upft ist) f\"uhrt zu den Isomorphismen
(der zweite f\"ur  $n$ prim zu Char($k$)): 
$$\br(k) \simeq H^2(\g,k_s^*)$$
$${}_n\br(k)=H^2(\g,\mu_n).$$ 
\bigskip

Die Brauergruppe tritt  in der Definition der
{\it Clifford Invariante von quadratischen Formen} vor (Varianten: Artin,
Hasse, Witt). 

%Zu jeder nichtentarteten quadratischen Form $q$ \"uber $k$
%der Dimension $d$ wird die Cliffordsche Algebra $C(q)$ assoziiert. 
%Die Dimension von $C(q)$ ist $2^d$. Man hat die gerade
%Cliffordsche Algebra $C_0(q)\quad\subset C(q)$, der Dimension
%$2^{d-1}$. Diejenige, deren  Dimension ein Quadrat ist, ist eine
%zentrale einfache Algebra. Ihre Klasse ist der (Clifford) Witt
%invariant $c(q) \in \br(k)$.

%Der Hassesche Invariant $s(q)$ 
%von 
%einer diagonalen Form $q=<a_1, \dots, a_n>$
%ist $\prod_{i< j} (a_i,a_j)$. Siehe Lam S. 123 f\"ur die Beziehung
%zwischen den beiden.

\bigskip

Sei $k$ algebraisch abgeschlossen und $K=k(C)$ der
Funktionen\-k\"orper einer Kurve. Dann ist $K$ ein $C_1$-K\"orper
und $\br(K)=0$ (Tsen).

\bigskip
%\vfill\eject

 Sei $A$ ein {\it diskret bewer\-teter henselscher Ring} (vom Rang Eins),
$K$ sei der Fraktionsk\"orper, $k$ der Residuenklassenk\"orper.
{\it Dann hat man eine exakte Folge
$$ 0 \to \br'(k) \to \br'(K) \to H^{1}{'}(k,\Q/\Z) \to 0.$$}
Strich hei{\ss}t: prim zur Charakteristik von $k$.
Diese Folge geht auf Witt zur\"uck. 

Sei $A$ nicht unbedingt henselsch. Elemente im Kern der
zusammengesetzten Abbildung
$\br'(K) \to \br'(\hat{K}) \to H^{1}{'}(k,\Q/\Z)$
nennt man {\it unverzweigt}.

\bigskip

 Die Brauergruppe von $k(t)$ 
($k$ perfekt) wird dann leicht berechnet (Faddeev). Man erh\"alt die
exakte Folgen
$$ 0 \to \br(k) \to \br(k(t)) \to \sum_{x \in \A^1_k{}^{(1)} }
H^1(k_x,\Q/\Z) \to 0$$
und
$$ 0 \to \br(k) \to \br(k(t)) \to \sum_{x \in \P^1_k{}^{(1)} }
 H^1(k_x,\Q/\Z) \to H^1(k,\Q/\Z) \to  0.$$

\bigskip

Unter Benutzung der Wittschen Folge und der Konstruktion von
Bewertungen auf Schief\-k\"orpern \"uber lokalen K\"orpern erh\"alt  man
die folgende sch\"one Formel (Tignol, ...).
Sei eine $K/k$ zyklische Erweiterung von K\"orpern,
und $A/k$ eine zentrale einfache Algebra. Dann ist 
$$\ind_{k(t) } (A \otimes_k(K/k,t))= [K:k]. \ind_K(A_K).$$

Anwendung ({\it Index kann gr\"osser als Exponent sein}):

 Wenn die
$K_i/k$ zyklisch und unabha\"ngig sind, und die $X_i$ Variablen, dann ist
die Algebra
$(K_1/k,X_1) \otimes \dots \otimes (K_n/k,X_n)$ \"uber dem
K\"orper $k(X_1,\dots,X_n)$ eine Divisionsalgebra (Nakayama, 1935;
spezielle F\"alle bei Brauer 1929 und bei K\"othe 1931).

%\bigskip
\vfill\eject

Eine {\bf Severi-Brauer Variet\"at} ist eine getwistete Form 
des projektiven Raumes (Beispiel: ein Kegelschnitt).
Solche Variet\"aten wurden bei beliebigen K\"orpern von Severi betrachtet,
bei Zahlk\"orpern von F. Ch\^atelet untersucht. Ch\^atelet nannte sie
``Brauer Variet\"aten". Der Name Severi-Brauer stammt aus einem
Bericht von B. Segre.

Sei $A/k$ eine zentrale einfache Algebra.
Die zugeh\"orige 
Severi-Brauer Variet\"at $X_A$ ist die
Variet\"at der (Rechts)idealen kleinster Dimension
von $A$. Wenn  $n$ der Index
von $A$ ist, so ist $X_A$ eine getwistete Form von $\P^{n-1}$.

Der Funktionenk\"orper $k(X_A)$ ist ein
 ``generischer Zerf\"allungsk\"orper" (Roquette, Amitsur).

\bigskip

Es gilt:
$$ Ker \hskip1mm [ \br(k) \to \br(k(X_A) ]= \Z [A].  $$

Allgemeiner, sei $X/k$ glatt und vollst\"andig, geometrisch
irre\-du\-zibel. Dann gibt es eine exakte Folge 
$$0 \to \pic(X) \to \pic({\overline X})^\g \to \br(k) \to \br(k(X)).$$

  Die klassischen S\"atze
\"uber Severi-Brauer Variet\"aten fliessen direkt aus
dieser exakten Folge.

% Qu'ont fait de plus Roquette et Amitsur ?

\bigskip

{\it Algebren mit Involution und klassissche Gruppen} (Weil)

\medskip

{\it Entwicklung der Galois-Kohomologie} (Serre, Tate)

\medskip

 {\it Wittgruppe} $W(k)$, Ideale $I^n(k)\subset W(k)$.

Die Arbeiten von Pfister und Arason \"uber quadratische Formen.
\hskip2mm Ker $W(k) \to W(k(X))$, wo $X/k$ eine Quadrik ist (Pfister).

\hskip2mm Hauptsatz von Arason-Pfister: $q$ anisotrop in
$I^nk$, dann ist dim($q$) wenigstens $2^n$.

\hskip2mm Invariante $e_3: I^3(k) \to H^3(k,\Z/2)$ (Arason).

Kern $H^3(k,\Z/2) \to H^3(k(X),\Z/2)$, f\"ur $X$ eine
Quadrik. (Arason)

\bigskip

{\it Nicht verschr\"ankte Produkte}: 
Amitsur (generisch); Saltman, Wadsworth,
Tignol, Jacob ; Brussel auf ${\bf Q}(t)$ und ${\bf Q}((t))$.

\bigskip

{\it ``Tannaka-Artin Problem"} (Platonov, Yanchevski\v{\i}, Draxl):
es gibt Beispiele von Algebren $A/k$, bei denen die Gruppe
$A^{*1}$ der Elemente von reduzierter Norm 1 gr\"o{\ss}er als die
Gruppe der Kommutatoren $[A^*,A^*]$ ist (der Quotient ist
als $SK_1(A)$ bekannt). Beispiele mit  $k=\Q_p((x))((y))$.

%Platonov und $SK_1$ ($Wh(G)$ ...) $SK_1(A)$, $Wh(G)$, $G(k)/R$.

%\vfill\eject
\bigskip

{\bf Unverzweigte Brauergruppe} 

\bigskip

Sei $K/k$ ein Funktionenk\"orper, Char.($k$)=0. Ein Element $\alpha$
von $\br(K)$ hei{\ss}t unverzweigt, wenn es bez\"uglich  jedem
diskret bewerteten (Rang Eins) Ring $A \subset K$ mit $k \subset A$
und Fraktionsk\"orper von $A$ gleich $K$ unverzweigt ist. Die Menge
aller solcher Elemente ist die unverzweigte Brauergruppe
$\br_{nr}(K/k)$ oder einfach $\br_{nr}(K)$ wenn der Grundk\"orper
$k$ klar ist.

Ist $K/k$ rein transzendent, so ist die nat\"urliche Abbildung
$\br(k) \to \br_{nr}(K/k)$ ein Isomorphismus.

Dies kann man benuzten, um rein Bewertungstheoretische 
Beweise der  Nichtrationalit\"at einiger unirationaler
Variet\"aten zu geben (Problem von L\"uroth; die Beispiele von
Artin-Mumford kann man so darstellen).
% CT-Ojanguren

\bigskip

 (Witt) {\it Sei $k=\R$, und $C/\R$ eine glatte zusammenh\"angende 
projektive Kurve. Dann ist
$$\br_{nr}(\R(C)) \simeq (\Z/2)^s$$
wo $s$ die Anzahl der Zusammenh\"angskomponenten
von  $C(\R)$  be\-zeich\-net.}

\vfill\eject
%\bigskip

(Bogomolov) {\it  Sei $L/{\bf C}$ ein Funktionenk\"orper beliebiger
Dimension. Sei $G$ eine endliche Gruppe von Automorphismen von $L$.
Dann ist
$$\br_{nr}(L^G)= \{  \alpha \in \br{L^G}, \alpha_H \in \br_{nr}(L^H)
\hskip1mm {\rm fuer }
%hskip1mm alle} 
\hskip1mm H \subset G \hskip1mm {\rm
bizyklisch}  \}.$$} 

\bigskip

(Bogomolov,
Saltman) {\it Sei $G \subset GL(V)$ eine endliche
Automorphismengruppe eines Vektorraumes $V/{\bf C}$. Dann
ist 
$$\br_{nr}(\C(V)^G)={ \rm Ker} \hskip1mm [ \hskip1mm H^3(G,\Z) \to
\prod_{A \subset G,
\hskip1mm A
\hskip1mm {\rm abelsch}} H^3(A,\Z) \hskip1mm ].$$} 

\bigskip

Auf diese Weise hat Saltman Beispiele von nichrationalen
 $\C(V)^G$ gegeben: der Noethersche Ansatz
zur Realisierung endlicher Gruppen als Galoisgruppen \"uber $\Q$
kann ohne weiteres nicht gelingen.

\bigskip

{\it Sei $G/\C$ eine reduktive zusammenh\"angende Gruppe,
und sei $G \subset GL(V)$ eine Einbettung. Dann
ist  $\br_{nr}(\C(V)^G)=0$.}

(Saltman f\"ur $G=PGL_n$; Bogomolov allgemein)

%\vfill\eject
\vskip1cm

{\bf II. Schematen.  Die Grothendiecksche Brauergruppe}

\bigskip

%Die Artikel von Grothendieck
%zum Thema Brauergruppe.

%\bigskip

Sei $X$ ein Schema. Dann kann man zwei Brauerartige Gruppen definieren:

\medskip

 Die {\it Azumaya(-Auslander-Goldman) Brauergruppe} 
$\br_{Az}(X)$: das ist die Gruppe der Klassen von Azumaya Algebren
auf $X$. Falls $A= {\cal {END}}  ({\cal E})$, mit ${\cal E}$ ein
Vektorb\"undel, dann ist $A$ trivial.

\medskip

{\it Grothendiecks Brauergruppe}: $\br(X)=H^2_{\et}(X,\G_m)$.

\medskip
Es gibt eine Einbettung $\br_{Az}(X) \hookrightarrow \br(X)$. Die erste
Gruppe ist eine Torsionsgruppe, die zweite nicht unbedingt.
Unter sehr allgemeinen Voraussetzungen hat Gabber angek\"undigt,
da{\ss} $\br_{Az}(X)$ mit der Torsion von $\br(X)$ \"ubereinstimmt.

Meistens ist man mit der Kohomologischen Gruppe $\br(X)$ ganz zufrieden,
besonders bei regul\"aren Schematen. 
%(F\"ur das Studium der
%$p$-Torsion kann manchmal der Azumaya Standpunkt von Hilfe  sein.)

\medskip

Lokalisierungsfolge  ($X$ regul\"ar)
$$0 \to \G_m \to i_{\eta}{_*}\G_m \to \oplus_{x \in X^{(1)}}i_{x*}{\bf Z}
\to 0.$$

Wenn $X$ regul\"ar und irreduzibel ist, mit Funktionenk\"orper $K$,
dann ist $\br(X) \to \br(K)$ injektiv (also ist $\br(X) $ torsion.)

% Dies ist speziell f\"ur $H^2(.,G_m)$.
% Dagegen erwartet ma solch eine Injektivit\"at
% im lokalen Fall f\"ur viele Funktoren (Gersten Vermutung)
% Clarifier le cas de $H^3(A,G_m) ...$
\medskip

F\"ur solch'ein $X$, darf man die Reinheitseigenschaft vermuten:
Sei $\xi \in \br(K)$. Nehmen wir an,  f\"ur alle $x \in X$ der
Kodimension Eins
 liegt $\xi$ im Bild von $\br(O_{X,x}) \to \br(K)$.
Geh\"ort $\xi$ zur Gruppe $\br(X) \subset \br(K)$~?
%% Noch mal Diskussion der Reinhaitsproblem, Gersten Vermutung
Dies ist in vielen F\"allen bewiesen, z.B. wenn $\dim(X) \leq 2$.
Wenn $X/k$ glatt ist und der K\"orper $k$ der Charakteristik null ist, so
gibt es eine exakte Folge $$ 0 \to \br(X) \to \br(K) \to \oplus_{x \in
X^{(1)}} H^1(k(x),\Q/\Z).$$
%F\"ur $k$ beliebig, stimmt das f\"ur die prim-zu-Char.($k$)-Torsion.

\bigskip

 Sei $A$ ein diskret bewer\-teter  (Rang Eins) Ring,
$K$ der Fraktionsk\"orper, $k$ der Residuen\-klassen\-k\"orper.
{\it Dann hat man eine exakte Folge} (Strich hei{\ss}t: prim zur
Charakteristik von $k$)
$$ 0 \to \br'(A) \to \br'(K) \to H^{1}{'}(k,\Q/\Z) \to 0.$$

\medskip

Bei glatten vollst\"andigen Variet\"aten $X/k$ (Char.($k$)=0)
folgt daraus die Birationale Invarianz von $\br(X)$.
$\br(X)$ stimmt mit der  unverzweigten Brauergruppe
von $k(X)$ (\"uber $k$)
\"uberein.

Man kann letztere benutzen, um bestimmte nicht triviale
Elemente vorzuf\"uhren. 

Die Grothendieck-Brauergruppe kann man dagegen
 global berechnen, wie wir in einer Minute  sehen werden.

\medskip

Grothendieck bietet h\"ohere birationale Invarianten an.
Sei $k$ ein K\"orper der Charakteristik null, und sei $X/k$ eine
irreduzibel vollst\"andige glatte Variet\"at.
F\"ur  $i,j \in \N$ h\"angt das Bild
der Restriktionsabbildung $$H^{i}(X,\mu_n^{\otimes j}) 
\to H^{i}(k(X),\mu_n^{\otimes j})$$
nur vom Funktionenk\"orper $k(X)$  ab. 
F\"ur $i=2, n=1$ erh\"alt man die $n$-torsion von $\br(X)$.

\bigskip

\bigskip

Sei $l$ eine Primzahl, $l$ eine Einheit auf $X$. Aus der Kummerschen
exakten Folgen :
$$1 \to \mu_{l^n} \to \G_m \to \G_m \to 1$$
($x \mapsto x^{l^n}$) erh\"alt man die exakte Folge
$$0 \to \pic(X)\otimes(\Q_l/\Z_l) \to
H^2_{\et}(X,\Q_l/\Z_l(1))
\to \br(X)\{l\} \to 0.$$

Falls $k=\C$, $X/\C$ projektiv glatt, dann ist 
$$ \br(X)\{l\} = (\Q_l/\Z_l)^{b_2-\rho} \oplus H^3(X(\C),\Z)_{{\rm
tors}}.$$

 $b_2-\rho >0$ ist zu 
$H^2(X,O_X)\neq 0$ \"aquivalent (Hodge-Theorie).

%\vfill\eject
\bigskip

Eine andere wichtige Methode zum Berechnen von Brauergruppen
liefert die Leray Spektralfolge. Sei $f: X \to Y$
ein Morphismus. Dann hat man eine Spektralfolge 
$$E_2^{pq}= H^p(Y,R^q\pi_*\G_m)
\Longrightarrow H^n(X,\G_m).$$
Nat\"urlich braucht man dann die Garben $R^q\pi_*\G_m$ 
zu berechnen.

\bigskip

Sei $X/k$ eine geometrisch
irreduzibel vollst\"andige glatte Variet\"at.  Sei $k_s$ ein
separabler Abschlu{\ss} von $k$, $\g=\gal(k_s/k)$,
$X_s=X\times_kk_s$. Dann hat man eine exakte Folge:
$$ 0 \to \pic(X) \to \pic(X_s)^\g \to \br(k) \to \hskip4cm $$
$$\hskip1cm \to 
Ker [\br(X) \to \br(X_s)] \to H^1(\g,\pic(X_s)) \to H^3(\g,k_s^*)$$

\medskip

Ist $X$ eine projektive Kurve
\"uber einem endlichen K\"orper $\F$, dann
folgt $\br(X)=0$.
Zusammen mit der Lokalisierungsfolge bekommt man daraus (bis auf die
$p$-Torsion) die exakte Folge der geometrischen Klassenk\"orpertheorie:
$$ 0 \to \br(\F(X)) \to \sum_{x \in X^{(1)}} \Q/\Z \to \Q/\Z \to 0.$$

\bigskip

 Sei $\pi: X \to Y$ eine eigentliche surjektive Abbildung,
sei $X$ regul\"ar, der Dimension 2,  sei $Y$ lokal henselsch.
Sei $Z$ die spezielle Faser. Unter angemessenen Annahmen,
$\br(X) \simeq \br(Z)$ (Satz von M. Artin: Verschwinden von
$R^2\pi_*\G_m)$.) 
 
\bigskip

Sei $\pi: X \to Y$ eine eigentliche surjektive Abbildung,
sei $X$ regul\"ar, der Dimension 2, sei $Y$ 
eine glatte vollst\"andige Kurve \"uber einem endlichen
K\"orper $\F$, und sei die generische Faser $X_{\eta}$ glatt und
geometrisch irreduzibel.  Unter
Benutzung des obigen Satzes wird die enge Beziehung zwischen
der Brauergruppe von $X$ und der Tate-Schafarewitsch-Gruppe 
der Jakobischen Variet\"at von $X_{\eta}$ festgestellt.

%\vfill\eject
\vskip1cm

{\bf III. $K$-Theorie (Grothendieck, Milnor, Bass, Quillen)}

\bigskip

Sei $X$ ein  Schema. Grothendieck hat
die Gruppe $K_0(X)$ der Klassen Vektorb\"undeln
auf $X$ definiert. F\"ur glatte Variet\"aten $X$ \"uber einem
K\"orper hat er die Beziehung zwischen $K_0(X)$ und den
Chowgruppen von $X$ untersucht 
(Chernklassen, Riemann-Rochscher Satz).

\bigskip

Sei $k$ ein K\"orper. Dann ist $K_0(k)=\Z$.
F\"ur $n \geq 1$ definiert man nach Milnor 
$$K_n^M(k) = (k^*\otimes_{\Z} \dots \otimes_{\Z} k^*)/ \{x_1\otimes \dots
\otimes x_n,
\hskip1mm x_i + x_j=1\}.$$
Also $K_1(k)=k^*$, 
und $K_2(k)=(k^*\otimes_{\Z} k^*)/ \{x\otimes y, x+y=1\}$.

\medskip

Es gibt nat\"urliche Homomorphismen ($p$ prim zu Char($k$))
$$K_n^M(k)/p \to H^n(k,\mu_p^{\otimes n})$$
 VERMUTUNG (Bloch-Kato) {\it  Diese Abbildungen  sind
Isomorphismen}
%(f\"ur $n=1$ hat man die Kummertheorie).

\medskip

Eine Folge  w\"are: die
nat\"urlichen Abbildungen $$H^{n+1}(k,\mu_r^{\otimes n})
\to H^{n+1}(k,\mu_{rs}^{\otimes n})$$ sind injektiv.

\medskip

Es gibt nat\"urliche Homomorphismen ($Char(k) \neq 2)$)
$$K_n^M(k)/2 \to I^n(k)/I^{n+1}(k).$$
 VERMUTUNG (Milnor) {\it Diese Abbildungen sind
Isomorphismen.}

%\vfill\eject
\bigskip
Quillen fand eine gute Definition f\"ur die h\"ohere
$K$-Theorie-Gruppen von Schematen. F\"ur diese Gruppen 
bewies er eine Reihe von Eingenschaften: Lokalisierung,
Homotopieinvarianz (im regul\"aren Fall), Gersten Vermutung
(bei glatten Variet\"aten),  Berechnung der h\"oheren
$K$-Theorie-Gruppen eines projektiven B\"undels, und auch 
eines Severi-Brauer Schemas (letzter Satz der Arbeit; der Fall $K_0$
ist schon interessant).

\bigskip

\"Ahnliche Berechnungen von $K_i(X,A)$ f\"ur $X/k$ projektiv homogen
und $A$ eine zentrale einfache  Algebra \"uber $k$ wurden sp\"ater
systematisch gemacht (Swan, Schofield-van den Bergh, Merkur'ev, Panin,
Wadsworth ...). Die Berechnnung von $K_0$ ist genau so schwer wie
die der h\"oheren $K_i$.

\bigskip

Man kann die $K$-Theorie eines Schemas mittels der Dimension
des Tr\"agers  koh\"arenter Garben filtrieren.
Es entsteht die Spektralfolge von Brown, Gersten, Quillen 
$$E_2^{pq}(X)= H^p_{Zar}(X,{\cal K}_q) \Longrightarrow K_n(X).$$

Auf die $K_i(X)$ kann man Chernklassen mit Werten in
der $K$-Kohomologie definieren, 
(Schechtman, Soul\'e, Kratzer)
und Gillet hat 
einen Riemann-Rochschen Satz bewiesen. Daraus bekommt man eine
partielle  Trivialisierung
% d\'eg\'en\'erescence
der BGQ Spektralfolge.

\bigskip
%\vfill\eject

{\bf Die S\"atze von Merkur'ev-Suslin}

\bigskip

Diese ganze Maschinerie haben Merkur'ev und Suslin benutzt.
In ihrer Arbeit wird die  
 Hauptrolle  durch die Brauer-Severi Variet\"at $X=SB(D)$
einer Divisionsalgebra $D$ vom Index eine Primzahl $p$
\"uber dem Grundk\"orper $k$ gespielt.

F\"ur diese Variet\"at $X$  beweisen Merkur'ev und Suslin:

\medskip

{\it $K_2k \to H^0(X,{\cal K}_2)$ ist ein Isomorphismus.

Die Restriktion
$H^1(X,{\cal K}_2) \to H^1(X_s,{\cal K}_2)$ ist eine Injektion.}

\medskip

Das sind sehr wichtige S\"atze auf dem Weg zu den Haupts\"atzen:

\medskip

{\it F\"ur $K/k$ zyklisch, $\gal(K/k)=\{\sigma\}$ und $K_2$
gilt ein Analog des Hilbertschen Satzes 90, n\"amlich 
$x \in K_2K$ kann man als $(\sigma-1)(y)$ schreiben, dann und
nur dann, wenn die Norm $N_{K/k}(x) \in K_2k$ verschwindet.}
\medskip

{\it Die Reziprozit\"atsabbildung 
$K_2^M(k)/p \to H^2(k,\mu_p^{\otimes 2})$
ist ein Isomorphismus.} 
\medskip

{\it Sei $k$ ein K\"orper, sei $A/k$ eine Algebra mit quadratfreiem Index
$n$. Dann liefert  die nat\"urliche Abbildung $f \mapsto f \cup [A]$
eine Einbettung $k^*/Nrd(A^*) \hookrightarrow H^3(k,\mu_n^{\otimes 2})$.}

\bigskip

{\bf Wie geht es weiter ?}

\bigskip

Die Geschichte l\"auft weiter (Rost, Merkur'ev-Suslin, Suslin, Voevodsky,
 Morel, Levine). Heutzutage wird die  $K$-Theorie durch die
motivische Kohomologie ersetzt (und die Brauergruppe $\br(k)$  schreibt
man jetzt
$H^3_L(k,\Z(1))$ oder $H^3_{\et}(k,\Z(1))$.)

\medskip

Im Prinzip ist die Bloch-Kato Vermutung f\"ur $p=2$ und $n$ beliebig
bewiesen, sowohl als auch der Isomorphismus $K_n^M(k)/2 \simeq
I^n(k)/I^{n+1}(k)$, also auch $I^n(k)/I^{n+1}(k) \simeq H^n(k,\Z/2))$.
So sind auch S\"atze, die die Ergebnisse von Arason verallgemeinern,
etwa: Kern von $H^n(k,\Z/2) \to H^n(k(X),\Z/2) =\Z/2$
wenn $X$ eine Quadrik ist, die durch eine anisotrope
$n$-Pfister Form definiert ist.

%\vfill\eject
\vskip1	cm
{\bf IV. H\"ohere Klassenk\"orpertheorie} (K. Kato, S. Saito)

\medskip

Aus der \"ublichen Klassenk\"orpertheorie hat man den Isomorphismus
 $$\br(k_v) \simeq H^1(\F_v,\Q/\Z) \simeq \Q/\Z$$
wo $k_v$ die Komplettierung eines globalen K\"orpers 
bzgl. einer diskreten Bewertung ist, und $\F_v$ den endlichen
Residuenklassenk\"orper bezeichnet. 
F\"ur einen globalen K\"orper $k$ hat man die exakte Folge
$$ 0 \to \br(k) \to \oplus_{v\in \Omega} \br(k_v) \to \Q/\Z \to 0.$$
Die Injektion $\br(k) \hookrightarrow \oplus_{v\in \Omega} \br(k_v)$
ist ein Haupsatz von Brauer-Hasse-Noether.
%``Jede \"uberall zerfallende
%Algebra u\"ber $\Omega$ ist $tilde ... \Omega$."
%% SIGNE tilde horizontal ...

\medskip

F\"ur Funktionenk\"orper von mehreren Variabeln \"uber einem
endlichen K\"orper oder einem Zahlk\"orper hat Kato 1986
Verallgemeinerungen dieser Tatsachen vermutet. Spezielle F\"alle  
sind bekannt (Kato, CT-Sansuc-Soul\'e, Jannsen, S. Saito, CT).

\vfill\eject

Spezielle Konsequenzen dieser Vermutungen:

\smallskip

A)  Sei $X/k$ eine geometrisch irreduzible Variet\"at der Dimension
$d$ \"uber einem Zahlk\"orper. Sei $k(X)$, bzw. $k_v(X)$, der
Funktionenk\"orper von $X$, bzw. $X \times_kk_v$.  Dann ist die
Abbildung
$$ H^{d+2}(k(X),\Q/\Z(d+1)) \to \prod_{v \in \Omega}
H^{d+2}(k_v(X),\Q/\Z(d+1))$$
eine Injektion.
F\"ur $d=0$ erkennt man den Satz von Hasse-Brauer-Noether.
Beweise f\"ur $d=1$ sind ver\"offentlicht; f\"ur $d$
beliebig seit einiger Zeit angek\"undigt.

\smallskip

B)  F\"ur $k$ lokal und $R$ den Ring der ganzen Elementen in $k$,
$\pi: {\cal X} \to {\rm Spec}(R)$ ein eigentlicher flacher  Morphismus,
${\cal X}$ regul\"ar, $X/k$ die generische Faser, der Dimension $d$, und
$Y/\F$ die spezielle Faser von $\pi$, Berechnung der $(d+2)$-unverzweigten
Kohomologie von $X/k$ (Koeffizienten $\Q_l/\Z_l(d+1)$) mittels der
$(d+1)$-Kohomologie von $Y/\F$ (Koeffizienten $\Q_l/\Z_l(d)$).
 Bekannt f\"ur $d=1$.

\vskip1cm

{\bf V. Ergebnisse}

\bigskip

\bigskip

{\bf Zerf\"allungsk\"orper}

\bigskip

%Algebren zerfallen nach einer aufl\"osbaren Erweiterung;
 {\it Eine Algebra $A/k$ vom Exponent $n$ (prim zu Char.($k$))
 besitzt
einen Zerf\"allungsk\"orper der Gestalt
$k(\mu_n)(a_1^{1/n},\dots,a_m^{1/n})$.} (Merkur'ev-Suslin)

\bigskip

{\bf Struktur der Brauergruppe eines K\"orpers} 

\medskip
{\it Falls $k$ alle Einheitswurzeln enth\"alt, so it
$\br(k)$ eine (Torsion) teilbare Gruppe.} (Merkur'ev-Suslin)

\medskip

{\it Die Gruppe ${}_p\br(k)$ wird durch Elemente vom Index $p$
erzeugt.} (Merkur'ev)

\medskip

{\it Sei $p\neq 2$. Wenn es eine nicht triviale zyklische
Algebra vom Index $p$ \"uber $k$ gibt, dann enth\"alt 
$\br(k)$ eine Untergruppe $\Q_p/\Z_p$.} (Merkur'ev)
% f\"ur p=2 \neq Char k, supposer que k(\mu_{E^{\infty}}/k$
% est sans torsion

\medskip

{\it Sei $p=2$ oder $3$. Wenn $2.\br(k)\{p\} \neq 0$, dann
enth\"alt $\br(k)$ eine Untergruppe $\Q_p/\Z_p$.}  (Merkur'ev)

%\vfill\eject
\bigskip
{\bf Index-Reduktion Formeln} 

\medskip

Die Methode, die von Schofield und van den Bergh stammt,
besteht darin, die Surjektivit\"at
$K_0(X,D) \to K_0(D_{k(X)}) \simeq \Z$ zu benutzen
($X/k$ glatt), und die Berechnung der Gruppen
$K_0(X,D)$ bei projektiven homogenen R\"aumen von
linearen algebraischen Gruppen auszuwerten.

\medskip

Urspr\"unglicher Satz (Schofield-van den Bergh):
{\it Seien $A$ und $B$ zwei
zentrale Algebren
\"uber
$k$. Sei $X_B$ die zu $B$ geh\"orige Severi-Brauer Variet\"at, und sei
$K_B=k(X_B)$. Es gilt:
 $$\ind (A\otimes_kK_B)= {\rm g.g.T.} \hskip1mm \ind(A\otimes_kB^{\otimes
n})= {\rm Inf} \hskip1mm \ind(A\otimes_kB^{\otimes n}).$$}

Erste Anwendungen (Schofield/van den Bergh):  Wildes verhalten des
Indexes bez\"uglich dem Exponent; Beispiele von  unzerlegbaren Algebren
mit kleinem Exponent.

Die Methode ist ganz anders als die Bewertungsmethoden, die
seit der Zeit von Brauer benutzt worden waren.
% geworden ... ?

\bigskip

{\bf $u$-Invariante von K\"orpern}

\medskip

 {\it  F\"ur jede gerade Zahl $2n$
gibt es einen K\"orper $k$, $cd(k) \leq 2$,  so da{\ss} das Maximum der
Dimensionen anisotroper quadratischer Formen \"uber $k$ gleich $2n$ ist.}
(Merkur'ev)

Beweis: Indexreduktionsformel und die Swansche Berechnung der Gruppen
$K_0(X,A)$ f\"ur $X$ eine Quadrik. Letztere kann man
durch 
einen elementaren Beweis (Tignol) ersetzen.

%\vfill\eject
\bigskip
{\bf Vermutung II von Serre}

\medskip

VERMUTUNG (Serre) {\it  F\"ur $G$ halbeinfach zusammenh\"angend  und
$k$ perfekt mit $cd(k)=2$,  ist $H^1(k,G)=0$.}
   
Bewiesen von Merkur'ev-Suslin f\"ur $G=SL(1,A)$
(f\"ur $A$ quadratfrei klar, siehe oben).
Bewiesen von Bayer-Parimala f\"ur  klassische Gruppen, $F_4$, $G_2$;
von Gille f\"ur quasi-zerfallende Gruppen ohne $E_8$.

\medskip

Variante der Vermutung II f\"ur K\"orper
mit $$cd(k(\sqrt{-1})) = 2$$ (lokal-globales Prinzip
f\"ur
$H^1(k,G)$  bzgl. den reellen Abschl\"ussen). 
F\"ur die klassischen Gruppen wurde diese Variante auch von 
Bayer-Parimala bewiesen.

\bigskip

{\bf Die Rostsche Invariante und andere Invarianten}

\medskip

Sei $G/k$ eine   geometrisch  fast einfache
Gruppe \"uber einem K\"orper $k$, und sei $G$
einfach zusammenh\"angend. Rost hat eine Invariante
$$H^1(k,G) \to H^3(k,\Q/\Z(2)).$$
definiert --  und studiert. 

Vorl\"aufer waren $H^1(k,SL(D))=k^*/Nrd(D^*) \to
H^3(k,\mu_n^{\otimes 2}) $ f\"ur $D/k$ zentrale einfache Algebra vom
Index
$n$, und die (viel schwerer zu definieren) Invariante $e_3$ von Arason
f\"ur $H^1(k,Spin(q))$.
 
% Quid de $H^1(k,SL(D)) \to H^3(k,_\mu_n^{\otimes 2}) $ ??
%En fait meme celui-la ne va pas sur un anneau $A$
% non local a cause de $H^1(A,GL(D))$.

%\medskip
% COMPLETER avec remarque de Serre
%F\"ur Gruppen vom Typ $F_4$ hat Rost andere Invarianten definiert,
%mit Werten in h\"oheren Galois Kohomologiegruppen $H^n(k,
%\bullet)$.

\medskip

Es sei hier bemerkt, da{\ss} diese Invarianten sich keineswegs
auf Invarianten von $H^1(A,G)$ f\"ur einen beliebigen (kommutativen)
Ring $A$ ausdehnen, ganz im Unterschied mit der Clifford Invariante 
der quadratischen Formen (die Invariante von Artin, Hasse und
Witt).

\bigskip
%\vfill\eject

{\bf Die Gruppe $SK_1(A)$}

\bigskip

{\it Es ist $SK_1(A)=0 $ wenn  $cd(k)=2$.} (Merkur'ev-Suslin).

\bigskip

Sei $A/k$ eine zentrale Divisionsalgebra,  $n$ sei ihr Index. 
Man kann einen nat\"urlichen Homomorphismus
$$ SK_1(A) \to H^4(k,\Q/\Z(3))/(K_2(k)\cup [A])$$
definieren (Suslin, Kahn).
Suslin verbindet solch'eine Abbildung mit 
den  Beispielen von Platonov ($SK_1(A)\neq 0$).

%Suslin vermutet, da{\ss} es einen nat\"urlichen Homomorphismus
%$$SK_1(A) \to
%H^4(k,\mu_n^{\otimes 3})/(H^2(k,\mu_n^{\otimes 2})  
%\cup [A])$$
%gibt, die bei bewerteten K\"orpern die Invariante
%von Platonov geben w\"urde. Er
% hat eine etwa schw\"achere Abbildung definiert.
%Bruno Kahn sagt mir, man kann eine nat\"urliche Abbildung

%\bigskip
\vfill\eject

Sei $A$ eine Biquaternionenalgebra, und sei $q$ eine
zugeh\"orige Albertform (das ist eine quadratische Form
der Dimension 6). Dann gibt es eine exakte Folge (Rost)
$$ 0 \to SK_1(A) \to H^4(k,\Z/2) \to H^4(k(q),\Z/2).$$

\medskip

Merkur'ev hat daraus eine \"ahnliche Darstellung f\"ur
$SK_1(A)$ bei beliebigen Algebren vom Index 4 hergeleitet,
im Besonderen eine Injektion

$$SK_1(A) \hookrightarrow H^4(k,\Z/2)/(H^2(k,\Z/2) \cup A^{\otimes 2}).$$

\bigskip

{\bf $R$-\"Aquivalenz und Rationalit\"at von algebraischen
Gruppen}

\medskip

Sei $G/k$ eine lineare algebraische Gruppe.
Die Punkte von $G(k)$, die man mittels einer
rationalen Kurve zu $1_G \in G(k)$ verbinden kann, bilden eine
normale Untergruppe $RG(k)\subset G(k)$. Schreiben wir
$G(k)/R=G(k)/RG(k)$. Ist $G$ eine $k$-rationale
Variet\"at, dann ist 
$G(K)/R$ f\"ur jede K\"orper\-erweiterung $K/k$.
Kurz gesagt, $G$ ist dann $R$-trivial.

\medskip

Aus der Eigenschaft $SK_1(A) \simeq SK_1(A\otimes_kk(t))$
(Bass, Platonov) ist Voskresenski\v{\i} zum folgenden Schlu{\ss}
gekommen: 
$$SK_1(A)= SL_{1,A}(k)/R.$$      

\medskip

Die Beispiele von $A/k$ mit $SK_1(A) \neq 0$ (Platonov) liefern also
Beispiele von Gruppen $G=SL_{1,A}/k$, die nicht $k$-rational sind.

\medskip

F\"ur jede Divisionsalgebra $A/k$ vom Index $4n$ hat Merkur'ev 
gezeigt: es gibt eine Erweiterung $K/k$ mit $SK_1(A_K) \neq 0$.
Also ist $SL_{1,A}$ nicht rational \"uber $k$ (schon
f\"ur $k$ ein Zahlk\"orper). Er benutzt die
obere exakte Folge von Rost und Indexreduktionsformeln,
sowie die Injektivit\"at von
$e_4$ auf Symbolen.

\bigskip

%Merkur'ev und Chernousov haben Formeln f\"ur 
%$G(k)/R$ f\"ur die andere einfach zusammenh\"angende 
%Gruppen vom Typ $A_n$
%berechnet.

%\medskip

Es ist relativ leicht, Beispiele von halbeinfachen
Gruppen zu geben, die weder einfach zusammenh\"angend
noch von adjunktem Typ sind, und die nicht rational sind
\"uber dem Grundk\"orper (die Nichtrationalit\"at kann man
mittels der unverzweigten Brauer Gruppe sehen).

\bigskip

Beispiele von nichtrationalen Gruppen von
adjunktem Typ wurden erst 1994 von Merkur'ev gegeben.
Daf\"ur berechnete er die Quotienten $G(k)/R$ 
der klassischen adjunkten Gruppen. Das benutzt relativ
elementare Methoden (auch von Gille entwickelt).

 Sei $q/k$ quadratische Form von gerader Dimension, 
 dann ist $$PSO(q)(k)/R= G_q(k)/k^{*2}.Hyp_q(k)$$
wo $G_q(k)$ die Gruppe der \"Ahnlichkeitsfaktoren von $q$, und
$Hyp_q(k) \subset k^*$ die Gruppe, die durch die $N_{E/k}(E^*)$
erzeugt ist, mit $q_E$ hyperbolisch.

Sei dim($q$) gerade, $d=disc_{\pm}(q) \in k^*$ kein Quadrat,
$Z=k(\sqrt{d})$.  Sei $C_0(q)/Z$ vom Index wenigstens 4. Dann ist
$PSO(q)$ nicht $k$-rational. 

Um dies zu zeigen wird ein Funktionenk\"orper
 $L/k$ gebaut, so
da{\ss} $PSO(q)(L)/R \neq 0$. Hier wird von Indexreduktionsformeln
Gebrauch gemacht.

Beispiele gibt es \"uber $k=\Q(t)$.

\vfill\eject
%\bigskip

{\bf H\"ohere unverzweigte Kohomologie}

\medskip

Sei $X/k$ eine glatte, irreduzible, projektive Variet\"at,
$k(X)$ sei ihr Funktionenk\"orper.
F\"ur $i\in \N$ und geeignete Koeffizienten $\mu$ kann man
die unverzweigte Kohomologiegruppe 
$$H^{i}_{nr}(k(X),\mu) \subset H^{i}(k(X),\mu)$$
genau so definieren wie die unverzweigte Brauergruppe von $k(X)$
(\"uber $k$). Diese Gruppe enth\"alt die h\"ohere birationale
Invariante von Grothendieck, ist aber f\"ur $i\geq 3$ meistens gr\"osser
(bei $k$ algebraisch abgeschlossen und $\mu$ endlich kann sie sogar
unendlich sein). Mit diesen Gruppen kann man die Nichtrationalit\"at
bestimmter  unirationaler Variet\"aten beweisen.
\medskip

{\bf H\"ohere unverzweigte Kohomologie
bestimmter $V/G$}

\medskip

F\"ur $G \subset GL(V)$ endlich, \"uber $\C$,
hat  Peyre vor kurzem  Beispiele
gefunden mit \break $H^2_{nr}(\C(V)^G,\Q/\Z)=0$ und
$H^3_{nr}(\C(V)^G,\Q/\Z)
\neq 0$. (Hier gibt es Beziehungen mit den Equivarianten Chow-Gruppen.)

\bigskip

Sei $G/k$ halbeinfach einfach
zusammenh\"angend,  $G\subset GL(V)$.
F\"ur solch'eine Gruppe $G$ sind
die Gruppen $H^1_{nr}(k(V/G),\Q/\Z)$ und
$H^2_{nr}(k(V/G),\Q/\Z(1))$ trivial.

Vor kurzem hat Merkur'ev  allgemeine Formeln f\"ur die Gruppe  
$H^3_{nr}(k(V/G), \Q/\Z(2))$  erhalten. Damit hat er
bei Gruppen vom Typ $A_n$ und vom Typ $D_n$ 
(f\"ur passende nicht algebraisch
abge\-schlos\-sene K\"orper $k$)
 Beispiele 
von  nichtrationalen $V/G$ gefunden. 
%(und Garibaldi)

% Saltman/Tignol f\"ur PGO ...)
\bigskip

{\bf Verallgemeinerung eines Satzes von Witt}

\medskip

Sei $X/\R$ eine glatte irreduzible projektive  Variet\"at
der Dimension $d$, und $s$ sei die Anzahl der Zusammenhangskomponenten
von $X(\R)$.
Dann ist 
$$H^{d+1}_{nr}(\R(X),\Z/2) \simeq (\Z/2)^s.$$
(CT-Parimala, weitere Arbeiten von Scheiderer)

%\vfill\eject
\bigskip

{\bf Arithmetik der zwei-dimensionalen henselschen Ringen}

\medskip

Sei $A$ ein 2-dimensionaler, lokaler, henselscher Integrit\"atbereich,
der  Residuenklassenk\"orper sei algebraisch abgeschlossen der Char.
null. Der  Fraktionsk\"orper $k$ benimmt sich fast genau wie ein total
imagin\"arer Zahlk\"orper: 
% expliquer les anneaux de valuation ...

{\it Alle Algebren sind zyklisch. F\"ur quadratische Formen
der Dimension 3 oder 4 gilt ein lokal-globales Prinzip; ist die Dimension
wenigstens 5, dann sind sie isotrop.} (Ford-Saltman;
CT-Ojanguren-Parimala)
% Vermutung II von Serre ...
% Maximale abelsche Erweiterung der cd 1 ..
%% HP f\"ur projektive homogene R\"aume ?

Hier werden  S\"atze von M. Artin benutzt; einer wurde schon
erw\"ahnt, der zweite (aus  SGA4) lautet:
$cd(k)=2$. 

\bigskip

{\bf Arithmetik der K\"orper $\Q_p(t)$}

\medskip

Algebren und quadratische Formen auf $k(C)$, mit $k$ 
lokaler K\"orper. Sagen wir $k$ $p$-adisch, $p \neq 2$.

\medskip

{\it Jedes Element von ${}_2\br(k(C))=H^2(k(C),\Z/2)$ 
zerf\"allt \"uber  einer Erweiterung der Gestalt
$k(C)(\sqrt{f},\sqrt{g})/k(C)$} (Saltman, benutzt Artin
 und die Tatsache da{\ss} die Brauergruppe  einer
voll\-kommenen Kurve \"uber einem endlichen K\"orper verschwindet).

\medskip

{\it Jedes Element von $H^3(k(C),\Z/2)$ zerf\"allt \"uber einer
quadratischen Erweiterung} (Parimala-Suresh; benutzt die h\"ohere
Klassen\-k\"orper\-theorie).

\medskip

{\it Jede quadratische Form in $n \geq 11$ Variabeln \"uber
$k(C)$ ist isotrop} (Parimala-Suresh, unter Benutzung der obigen
Ergebnisse, sowohl als auch F\"alle der Bloch-Kato Vermutung; vorher
Merkur'ev, Hoffmann-van Geel).
% VERIFIER

%? Conjecture de Gersten (injection, puret\'e) 
%pour $H^1(., SL(D))$ ? Ct-Oj, 
% CT-Parimala-Sridharan, Rost, Panin-Suslin

%\vfill\eject

\vskip1cm
{\bf VI. Einige offene Probleme}

\bigskip

{\bf Zerf\"allungsk\"orper, Index und Exponent}

\medskip

{\it Wenn  $\ind(A)$ eine Primzahl $p$ ist,  ist $A$ zyklisch ?}

Bekannt f\"ur $p=2,3$. Unbekannt f\"ur $p=5$ (sogar wenn  $k={\bf
C}(x,y)$).

\bigskip

{\it Wenn $k$ ein $C_2$-K\"orper ist, ist $\ind(A)=\exp(A)$ ?}

Der Fall, in dem $\exp(A)$  eine Potenz von 6 teilt, ist bekannt (Artin, 
Harris, Tate, Yanchevski\v{\i}, Merkur'ev-Suslin).
F\"ur  beliebige Exponente ist die Frage schon offen wenn $k={\bf
C}(x,y)$.

%Quid de corps interm\'ediaires
%type ${\bf C}((x))(y)$ ou ${\bf C}(x)((y))$ ?

\bigskip

 Sei $k={\bf C}(x_1,\dots,x_d)$. Sei $p$ eine Primzahl,
$\zeta_p$ eine primitive Einheitswurzel. Das Tensorprodukt
$$(x_1+2,x_2)\otimes(x_1+3,x_3)\otimes \dots \otimes (x_1+d,x_d)$$
ist eine Divisionsalgebra vom Exponent $p$ und vom Index $p^{d-1}$.

{\it  Kann man Divisionsalgebren \"uber $k$ vom Exponent $p$ und von
gr\"osserem Index finden ?}

\bigskip

{\it Ein Beispiel wo 
 $k$ ein Funktionenk\"orper vom Transzendenzgrad $d$ 
\"uber ${\bf C}$ ist,  $A/k$ unverzweigt und vom Exponent $p$ ist, und
$\ind(A) \geq p^{d-1}$  geben.}

%Kandidat: der Funktionenk\"orper von $X=C\times E_1 \times \dots \times
%E_{d-1},$ wo  $C$ eine Kurve von
%hohem Geschlecht und die $E_i$ passende elliptische Kurven sind.

\bigskip

{\bf   Algebren \"uber $k(t)$} 

\medskip

{\it Sei $k$ ein K\"orper und $t$ eine Variabel.
Sei $A \in \br(k(t))$, sei
$A$ in $t=0$ unverzweigt und $A(0)=0 \in \br(k)$.
  Gibt es $t=f(s) \in k(s)$, wenn m\"oglich mit $f(0)=0$,
so da{\ss} $A \otimes_{k(t)}k(s)=0 \in \br(k(s))$  ?}

\medskip

Iskovskih, Mestre: Ja wenn $A$ 2-torsion und wenig verzweigt ist.

Yanchevski\v{\i}:  Ja wenn $k$ lokal (oder sogar ``large") ist.
% Lien avec $G$ groupe de Galois sur Q ?!!!
\bigskip
%\vfill\eject

{\bf Die Gruppe $SK_1(A)$}

\medskip 

Sei $A/k$ eine zentrale einfache Algebra.

\bigskip

(Suslin)  {\it Wenn $cd(k)=3$, ist $SK_1(A)=0$ ?} 

Ja wenn $\ind(A)=4 $ (Rost, Merkur'ev)

\bigskip

{\it Sei  $k$ endlich erzeugt \"uber dem Grundk\"orper. 
Ist $SK_1(A)$ endlich ?}

Ja wenn $A$ eine Biquaternionenalgebra ist, und $k$ ist vom
Transzendenzgrad 2 \"uber einem globalen K\"orper.

\bigskip

{\bf Andere lineare Gruppen}

\medskip

Indem $SK_1(A)=SL(1,A)/R$, liegt es nahe, \"ahnliche
Probleme f\"ur $G(k)/R$ bei anderen zusammenh\"angenden linearen
algebraischen Gruppen
$G$ aufzuwerfen. 
\medskip

{\it Ist $G(k)/R$  abelsch ?} 

%Trivial wenn $G/k$ zusammenh\"angend und $c(d(k)=1$.
% $G/k$ einfach zusammenh\"angend  und  $cd(k)=2$: ist $G(k)/R=1$ ?

\medskip

{\it Sei $G/k$ einfach zusammenh\"angend  und  $cd(k) \leq 3$. Ist dann 
$G(k)/R=1$ ? }
% Ja wenn $cd(k)=1$. Schon klar bei jedem zusammenh. Gruppe ..

\medskip

{\it Sei $k$ endlich erzeugt \"uber dem Primk\"orper. Ist $G(k)/R$ endlich
?}

Bekannt f\"ur $k$ ein Zahlk\"orper (Gille).

%\bigskip
\vfill\eject

{\bf  Normenprinzip}

\medskip

 Sei $A/k$ zentral einfach.
Sei $K/k$ eine endliche Erweiterung.
Dann ist $N_{K/k}(Nred(A_K^*)) \subset Nred(A^*) \subset k^*$.
Bei quadratischen Formen haben Knebusch und Scharlau andere
Normenprinzipen bewiesen.

Man kann fragen, ob es sich um ein allgemeines Ph\"anomen in der Theorie
der  algebraischen Gruppen handelt.

 Betrachten wir exakte Folgen 
$$ 1 \to \mu \to \tilde{G} \to G \to 1$$
$$ 1 \to H \to G \to T \to 1 $$
wo $\mu$ endlich von multiplikativem Typ ist,
$\tilde{G}, G$ zusammen\-h\"angende reduktive 
Gruppen sind, und $T$ ein Torus ist.

Gille und (sp\"ater) Merkur'ev haben die Fragen untersucht {\it ob ein
Normenprinzip f\"ur das Bild von $G(k) \to H^1(k,\mu)$ und f\"ur das Bild
von $G(k) \to T(k)$  gilt.} 

Gille und Merkur'ev zeigen, da{\ss} ein Normenprinzip gilt, wenn man
sich  auf  das Bild der Untergruppe $RG(k) \subset G(k)$ beschr\"ankt.
%dabei bezeichet
%$RG(k)$  die (normale) Untergruppe der Elementen von $G(k)$, die
%$R$-\"aquivalent zu $1 \in G(k)$ sind. 
Also wenn $G$ $k$-rational ist,
dann gilt das Prinzip. Dies liefert einen neuen Beweis des Normenprinzips
von Knebusch,  f\"ur das Normenprinzip von Scharlau reicht die Methode
aber nicht aus.

% Merkur'ev-Chernousov: Normenprinzip f\"ur G(k)/R,
% klassische Gruppen
%% Exemples Merkur'ev pour $\lambda: G \to T$.
%% d\'etailler

\bigskip

{\bf Rationalit\"at des Zentrums der generischen Algebra}

\medskip

{\it Sei $V/\C$ eine treue lineare Darstellung von $PGL_{n,\C}$.
Ist der Quotient $V/PGL_n$ eine rationale Variet\"at ?}

\medskip

$H^2$ unverzweigt ist null (Saltman).
% autre d\'emonstration CT-Sansuc
 Dasselbe gilt f\"ur $H^3$ unverzweigt (Saltman).
% VERIFIER 

\medskip

Man kennt sowieso kein Beispiel einer zusammenh\"angenden reduktiven
Gruppe
$G \subset GL(V)$
\"uber $\C$ mit $V/G$ nicht rational.

%\vfill\eject
\bigskip

{\bf Unverzweigte Kohomologie}

\bigskip

FRAGE {\it Sei $X/\F$ eine irreduzible glatte Variet\"at
\"uber einem endlichen K\"orper $\F$. F\"ur $i,j,n \in \N$,
mit $n$ prim zur Charakteristik von $\F$, ist die 
unverzweigte Kohomologiegruppe
$$H^{i}_{nr}(\F(X),\mu_n^{\otimes j})$$
endlich ?}

\bigskip

FRAGE {\it Sei $X/k$ eine irreduzible glatte Variet\"at
\"uber einem Zahl\"orper $k$. F\"ur $v$ Stelle von $k$,
sei $k_v(X)$ der Funktionenk\"orper von $X_v=X\times_kk_v$.
F\"ur $i,j,n \in \N$, ist der Kern der diagonalen Abbildung
$$H^{i}(k(X),\mu_n^{\otimes j}) \to \prod_v
H^{i}(k_v(X),\mu_n^{\otimes j}) $$
endlich ?}

\medskip

F\"ur $i\geq 3$ scheint es gewagt, die selbe Frage mit 
Koeffizienten $\Q/\Z(j)$ aufzuwerfen. 

F\"ur $i=2$ ist das aber im Grunde genommen
eine wohlbekannte Frage.

\vfill\eject

{\bf Eine Vermutung von Tate}

\medskip

Im Rahmen seiner allgemeinen Vermutungen \"uber algebraische
Zyklen hat Tate die folgende Vermutung aufgeworfen:

\bigskip

VERMUTUNG (Tate) {\it Sei $X/\F$ eine glatte vollst\"andige
Variet\"at \"uber einem endlichen K\"orper. 
Dann ist  $\br(X)$ 
endlich.}

\bigskip

In dem Fall, wo $X$ eine Fl\"ache ist, mit einer Faserung $\pi: X \to
C$ \"uber einer Kurve, die Endlichkeit von $\br(X)$ ist 
eng mit der Endlichkeit der
Tate-Schafarewitsch Gruppe der generischen Faser von $\pi$ verkn\"upft.

\bigskip

Spezielle F\"alle sind bekannt:
 Rationale Fl\"achen (trivial),
 Abel\-sche Variet\"aten und 
 Produkte von Kurven (Tate),
spezielle $K3$-Fl\"achen (Artin und
Swin\-ner\-ton-Dyer).

\vskip1cm
%\vfill\eject

{\bf VII. Das Brauer-Manin-Hindernis zum Hasse-Prin\-zip
f\"ur rationale Punkte}

\bigskip

Sei $k$ ein Zahlk\"orper, $\Omega$ die Menge seiner Stellen.
Sei $X/k$ eine projektive, glatte, geometrisch irreduzible Variet\"at
\"uber $k$.

Man kann die Menge $X({\bf A}_k)$ der adelischen Punkte von $X$
mit der Brauergruppe paaren:

$$ X({\bf A}_k) \times \br(X) \to \Q/\Z .$$

 Das Bild der diagonalen Einbettung $X(k)
\to X({\bf A}_k)$ liegt im Kern $X({\bf A}_k)^{\br}$ der obigen Paarung
(Manin, 1970).

(Dabei wird von der exakten Folge der Klassenk\"orpertheorie
$$  0 \to \br(k) \to \oplus \br(k_v) \to \Q/\Z \to 0$$
 Gebrauch gemacht.)

% mention des fonctions zetas et du resultat de Brauer
% chez Cassels-Guy pour leur contre-exemple ?

\medskip

FRAGE \hskip2mm {\it Falls $X({\bf A}_k)^{\br} \neq \emptyset$,
ist auch $X(k) \neq \emptyset$ ?}

\bigskip

F\"ur Kurven vom Geschlecht Eins, z.B.
$ y^2=P(t)$ mit  $P(t)$
vom Grad 4,  folgt dies
aus der  Annahme, da{\ss} die Tate-Schafarewitsch
Gruppe (der Jacobischen Kurve) endlich ist. \"Ahnliches gilt f\"ur
prinzipielle homogene R\"aume von abelschen Variet\"aten.

\medskip

F\"ur  homogene R\"aume zusammenh\"angender
linearer algebraischen Gruppen, mit zusammenh\"angenden
(geometrischen) Isotropie\-gruppen,
hat die Frage eine bejahende Antwort
(Hasse, Schilling, Maa{\ss},
% lui aussi ??
 Eichler, Landherr, Kneser,
Harder, Tchernousov, Sansuc, Borovoi).

\medskip

F\"ur Fl\"achen der Gestalt $y^2-az^2=P(t)$ mit $P(t)$
vom Grad 4 wurde die Frage 1984 bejaht (CT, Sansuc und
Swin\-ner\-ton-Dyer).

%\vfill\eject
\bigskip

F\"ur beliebige (projektive, glatte) Variet\"aten war eine positive
Antwort {\it nicht  zu  erwarten}. Man musste
% mu{\ss}te ?
 allerdings bis 1998 warten, bevor ein explizites
(Gegen-)Beispiel gegeben wurde 
(Skorobogatov; weitere Arbeiten von Harari
und Skorobogatov). Bei diesen Beispielen ist die geo\-me\-trische
Fundamentalgruppe nicht trivial (sogar nicht abelsch).  

\bigskip
\vfill\eject

VERMUTUNG {\it Sei $X$ eine (glatte, vollst\"andige) Variet\"at.
Neh\-men wir an,  es gibt einen surjektiven Morphismus $p: X \to {\bf
P}^1_k$,  und folgende Bedingungen sind erf\"ullt~:

(a) die allgemeine Faser von $p$ ist birational zu einem prinzi\-piellen
homogenen Raum unter einer zusammen\-h\"an\-genden algebraischen Gruppe
$G$ \"uber dem K\"orper $k(\P^1)$, mit zusammen\-h\"an\-genden
(geometrischen) Isotropiegruppen;

(b)  f\"ur irgendeinen geschlossenen Punkt $M \in \P^1_k$ besitzt
die Faser $X_M$ eine Komponente der Multiplizit\"at Eins.

Wenn  $X({\bf A}_k)^{\br} \neq \emptyset$,
dann ist auch $X(k) \neq \emptyset$.}

\bigskip

F\"ur Familien $X \to {\bf P}^1_k$ von Kegelschnitten, also
 Gleichungen 
$$ a(\lambda)X^2+b(\lambda)Y^2+c(\lambda)T^2=0$$
gibt es  gute Gr\"unde, so was zu vermuten. 

Solch'eine Fl\"ache ist eine (geometrisch) rationale
Fl\"ache. F\"ur eine rationale Fl\"ache $X$ vermutet man sogar,
da{\ss} die Menge $X(k)$ 
dicht in  $X({\bf A}_k)^{\br}$ liegt.

%Allgemeine Gleichungen der Gestalt
%$$Norm_{K/k}(x_1\omega_1+ \dots + x_n\omega_n)=P(t),$$
%wo $K/k$ eine endliche aber nicht unbedingt zyklische
%Erweiterung von $k$ vom Grad $n$ sind noch r\"atselhaft.

%\vfill\eject
\bigskip

Vor kurzem wurden spezielle Fl\"achen, die eine Schar
% voir chez Max Noether
 von Kurven vom
Gesch\-lecht Eins besitzen, untersucht (Swin\-ner\-ton-Dyer, CT und
Skorobogatov). Eine subtile Variation hat Swin\-ner\-ton-Dyer (2000) bei
diagonalen kubischen Fl\"achen (das sind rationale Fl\"achen) entwickelt.
Dies f\"uhrt zum Satz:

\medskip

(Swinnerton-Dyer) {\it Nehmen wir an, die Tate-Schafarewitsch Gruppen
von elliptischen Kurven \"uber einem Zahlk\"orper sind
% seien ?
endlich. Dann gilt das Hasse-Prinzip f\"ur diagonale kubische
Hyperfl\"achen
$\sum_{i=0}^4 a_ix_i^3=0$
\"uber $\Q$.}

\bigskip

Ersetzt man $\bf Q$ durch den Funktionenk\"orper $\F_p(C)$
einer Kurve \"uber $\F_p$ mit $p\cong 2
\hskip1mm {\rm mod } \hskip1mm 3$, dann erh\"ahlt man
einen absoluten Satz: denn die Annahme bzgl. den
Tate-Schafarewitsch Gruppen wird  durch bekannte
F\"alle der
Tateschen Vermutung ersetzt.

%\vfill\eject
\bigskip

{\bf Das Brauer-Manin-Hindernis   zum
Hasse-Prin\-zip f\"ur Nullzyklen}

\bigskip

Sei $k$ ein K\"orper und $X$ eine $k$-Variet\"at.
Ein Nullzyklus von $X/k$ ist eine ganzahlige Kombination $\sum_Pn_PP$
von geschlossenen Punkten. Der Grad des Nullzyklus 
$\sum_Pn_PP$ ist die ganze Zahl $\sum_P n_P [k(P):k]$.
Ein Element $A \in \br(X)$ kann man auf $\sum_Pn_PP$
aufwerten:
$$ <A, \sum_P n_P P> = \sum_P n_P {\rm Cores}_{k(P)/k}(A(P)) \in \br(k)$$
(hier liegt $A(P) \in \br(k(P))$.)

\bigskip

Sei $k$ ein Zahlk\"orper, $
\Omega$ die Menge seiner Stellen.

\medskip

VERMUTUNG {\it 
Sei $X/k$ eine projektive, glatte,
geometrisch irreduzible Variet\"at \"uber $k$.
Nehmen wir an, auf $X$  gibt es eine Familie $\{z_v\}_{v \in \Omega}$
von  Nullzyklen vom Grad Eins, mit $z_v$ Nullzyklus auf $X
\times_kk_v$,
 so da{\ss} 
$ \sum_v <A,z_v>= 0$ f\"ur alle $ A \in
\br(X).$
Dann gibt es auf $X$ ein Nullzyklus vom
Grad Eins.}

\bigskip

F\"ur $X/k$ eine Kegelschnittfamilie \"uber $\P^1_k$ wurde die 
Vermutung 1988 bewiesen (Salberger). Der Beweis dehnt sich 
zu Familien von Severi-Brauer Variet\"aten \"uber $\P^1_k$ aus.

F\"ur $X$ eine Kurve von beliebiegem Geschlecht folgt die Vermutung
aus der Annahme, da{\ss} die Tate-Schafarewitsch
Gruppen von Abelschen Variet\"aten endlich sind (Manin, S. Saito).

F\"ur $X$ eine Kegelschnittfamilie \"uber einer Kurve von beliebiegem
Geschlecht folgt die Vermutung
aus der  selben Annahme (CT 2000; Frossard 2001).

\bye